# On the Basis Property of the Root Functions of Differential Operators with Matrix Coefficients

O. A. Veliev

Department of Math., Faculty of Arts and Sci., Dogus University,

Acıbadem, 34722, Kadıköy, Istanbul, Turkey.

e-mail: oveliev@dogus.edu.tr

**Abstract**

We obtain asymptotic formulas for eigenvalues and eigenfunctions of the operator generated by a system of ordinary differential equations with summable coefficients and periodic or antiperiodic boundary conditions. Then using these asymptotic formulas, we find necessary and sufficient conditions on the coefficients for which the system of eigenfunctions and associated functions of the operator under consideration forms a Riesz basis.

## 1 Introduction

Let $L(P_2, P_3, ..., P_n)$ be the operator generated in $L_2^m[0,1]$ by the differential expression

$$l(y) = y^{(n)}(x) + P_2(x) y^{(n-2)}(x) + P_3(x) y^{(n-3)}(x) + ... + P_n(x)y(x) \qquad (1)$$

and the periodic boundary conditions

$$y^{(\nu)}(1) = y^{(\nu)}(0), \ \nu = 0, 1, ..., (n-1), \qquad (2)$$

where $n$ is an even integer, $P_\nu(x) = (p_{\nu,i,j}(x))$ is a $m \times m$ matrix with the complex-valued summable entries $p_{\nu,i,j}(x)$ for $\nu = 2, 3, ...n$. Here $L_2^m[0,1]$ is the space of the vector functions $f = (f_1, f_2, ..., f_m)$, where $f_k \in L_2[0,1]$ for $k = 1, 2, ..., m$. In $L_2^m[0,1]$ the norm $\|.\|$ and inner product $(.,.)$ is defined by

$$\|f\|^2 = \int_0^1 |f(x)|^2 \, dx, \ (f,g) = \int_0^1 \langle f(x), g(x) \rangle \, dx,$$

where $|.|$ and $\langle .,. \rangle$ are the norm and inner product in $\mathbb{C}^m$. We often write $L$ for $L(P_2, P_3, ..., P_n)$. In this paper we obtain asymptotic formulas for the eigenvalues and eigenfunctions and then find necessary and sufficient conditions on the coefficient $P_2(x)$ for which the system of the eigenfunctions and the associated functions (root functions) of the operator $L$ forms a Riesz basis in $L_2^m[0,1]$. We shall work only with the periodic problem (1), (2). The changes which have to be done for the antiperiodic problem are obvious, and we shall note on them at the end of the paper.

First we discuss the papers devoted to the basis property of the root functions of the Sturm-Liouville operator $H$ generated in $L_2[0,1]$ by the differential expression

$-y''(x) + q(x)y(x)$ and the periodic boundary conditions, i.e., we discuss the case $n = 2$, $m = 1$. For brevity, we discuss only the periodic problem. The antiperiodic problem is similar to the periodic problem. It is known [1, Chap. 2] that the operator $H$ is regular but not strongly regular. The root functions of the strongly regular differential operators form a Riesz basis (this result is proved independently in [2-4]). In the case when an operator is regular but not strongly regular the root functions, generally, do not form even usual basis. However, it is known [5,6] that they can be combined in pairs, so that the corresponding 2-dimensional subspaces form a Riesz basis of subspaces (see for the definitions of the Riesz basis of subspaces in [7, Chap. 6], for example). In 1996 at the seminar in MSU Shkalikov formulated the following result. Assume that $q(x)$ is a smooth potential, $q^{(k)}(0) = q^{(k)}(1) = 0$ for $0 \leq k \leq s-1$, and $q^{(s)}(0) \neq q^{(s)}(1)$. Then the root functions of the operator $H$ form a Riesz basis in $L_2[0,1]$. Kerimov and Mamedov [8] obtained the rigorous proof of this result in the case $q \in C^4[0,1]$, $q(1) \neq q(0)$. Actually, this results remains valid for an arbitrary $s \geq 0$. It is obtained in Corollary 2 of [9].

Another approach is due to Dernek and Veliev [10]. The result was obtained in terms of the Fourier coefficients of the potential $q$. Namely, we proved that if conditions

$$\lim_{n \to \infty} \frac{\ln |n|}{n q_{2n}} = 0, \tag{3}$$

$$q_{2n} \sim q_{-2n} \tag{4}$$

hold, then the root functions of $H$ form a Riesz basis in $L_2[0,1]$, where $q_n = (q, e^{2\pi i n x})$ is the Fourier coefficients of $q$ and $a_n \sim b_n$ means that $a_n = O(b_n)$ and $b_n = O(a_n)$ as $n \to \infty$. Makin [11] improved this result. Using another method he proved that the assertion on the Riesz basis property remains valid if condition (4) holds, but condition (3) is replaced by a less restrictive one: $q \in W_1^s[0,1]$,

$$q^{(k)}(0) = q^{(k)}(1), \quad \forall k = 0, 1, ..., s-1, \tag{5}$$

$$|q_{2n}| > c_0 n^{-s-1}, \quad \forall n \gg 1 \quad \text{with some } c_0 > 0,$$

where $s$ is a nonnegative integer. Moreover, some conditions which imply the absence of the Riesz basis property were presented in [11]. Some sharp results on the absence of the Riesz basis property were obtained by Djakov and Mitjagin [12].

The results which we obtained in [9] are more general and cover all the previous ones except constructions in [12]. Several theorems on the Riesz basis property of the root functions of the operator $H$ are proved. One of the main results of [9] is the following. Let $q$ belong to the Sobolev space $W_1^p[0,1]$ with some integer $p \geq 0$ and satisfy condition (5), where $s \leq p$. Let the functions $Q$ and $S$ be defined by the equalities

$$Q(x) = \int_0^x q(t)\,dt, \ S(x) = Q^2(x)$$

and let $q_n, Q_n, S_n$ be the Fourier coefficients of $q, Q, S$ with respect to the trigonometric system $\{e^{2\pi i n x}\}_{-\infty}^{\infty}$. Assume that the sequence $q_{2n} - S_{2n} + 2Q_0 Q_{2n}$ decreases not faster than the powers $n^{-s-2}$. Then the root functions of the operator $H$ form a Riesz basis in the space $L_2[0,1]$ if and only if the following condition holds

$$q_{2n} - S_{2n} + Q_0 Q_{2n} \sim q_{-2n} - S_{-2n} + 2Q_0 Q_{-2n}.$$

If $n = 2\mu + 1$ and $m = 1$, then the operator $L$ is strongly regular and hence its root functions form a Riesz basis (see [1-4] ). The case $n = 2\mu + 1 > 1$ and $m$ is an arbitrary integer is investigated in [13], where we proved that if the eigenvalues $\mu_1, \mu_2, ..., \mu_m$ of the matrix

$$C = \int_0^1 P_2(x)\,dx \tag{6}$$

are simple, then the eigenvalues of $L$ are asymptotically simple and the root functions form a Riesz basis.

In this paper we consider the case $n = 2\mu$ and $m$ is an arbitrary integer. This case is more complicated, since even in the simple subcase $m = 1$ operator $L$ is not strongly regular. Moreover, the simplicity of the eigenvalues $\mu_1, \mu_2, ..., \mu_m$ does not imply that the eigenvalues of $L$ are asymptotically simple. First we obtain asymptotic formulas for the eigenvalues and eigenfunctions of $L$. Then we find necessary and sufficient conditions on the coefficient $P_2(x)$ for which the root functions of the operator $L$ form a Riesz basis in $L_2^m[0,1]$. To describe the conditions on $P_2(x)$ let us introduce some notations. Let $v_1, v_2, ..., v_m$ be the normalized eigenvectors of the matrix $C$ corresponding to the eigenvalues $\mu_1, \mu_2, ..., \mu_m$. Denote by $w_j$ for $j = 1, 2, ..., m$ the eigenvector of the adjoint matrix $C^*$ corresponding to $\overline{\mu_j}$ and satisfying $(w_j, v_j) = 1$. Introduce the notations

$$b_{s,q}(x) = \langle P_2(x)v_q, w_s\rangle, \ b_{s,q,p} = \int_0^1 b_{s,q}(x) e^{-2\pi i p x} dx, \ b_k = \max_{i,j=1,2,...,m}\{|b_{i,j,k}|\}. \quad (7)$$

In this paper we prove that if

$$\lim_{k\to\infty} \frac{\ln|k|}{k b_{s,s,\pm 2k}} = \lim_{k\to\infty} \frac{b_{2k}b_{-2k}}{b_{s,s,\pm 2k}} = 0, \ \forall s, \quad (8)$$

then the root functions of $L$ form a Riesz basis if and only if $b_{s,s,2k} \sim b_{s,s,-2k}$ for all $s = 1, 2, ..., m$. The similar results are obtained for the operator $A$ generated by (1) and the antiperiodic boundary conditions

$$y^{(\nu)}(1) = -y^{(\nu)}(0), \ \nu = 0, 1, ..., (n-1). \quad (9)$$

Let us introduce some preliminary results and describe the scheme of the paper. Clearly,

$$e_j e^{\pm i 2\pi k x} \text{ for } j = 1, 2, ..., m, \text{ where } e_1 = \begin{pmatrix} 1 \\ 0 \\ \vdots \\ 0 \end{pmatrix}, e_2 = \begin{pmatrix} 0 \\ 1 \\ \vdots \\ 0 \end{pmatrix}, ..., e_m = \begin{pmatrix} 0 \\ \vdots \\ 0 \\ 1 \end{pmatrix}$$

and $k \in \mathbb{Z}$, are the normalized eigenfunctions of the operator $L(0)$ corresponding to the eigenvalue $(2\pi k i)^n$. Here the operator $L(P_2, P_3, ..., P_n)$ is denoted by $L(0)$ when $P_v(x) = 0$ for $v = 2, 3, ..., n$. It easily follows from the classical investigations [1, chapter 3, theorem 2] that boundary condition (2) is regular and all large eigenvalues of $L$ consist of the sequences

$$\{\lambda_{k,1} :| k |\geq N\}, \{\lambda_{k,2} :| k |\geq N\}, ..., \{\lambda_{k,m} :| k |\geq N\}, \quad (10)$$

where $N \gg 1$, $k \in \mathbb{Z}$, satisfying the following asymptotic formulas

$$\lambda_{k,j} = (2\pi k i)^n + O\left(k^{n-1-\frac{1}{2m}}\right), \ \forall j = 1, 2, ..., m. \quad (11)$$

The method proposed here allows us to obtain the asymptotic formulas of high accuracy for the eigenvalue $\lambda_{k,j}$ and for the corresponding normalized eigenfunction $\Psi_{k,j}(x)$ of $L$ when $p_{\nu,i,s} \in L_1[0,1]$ for all $\nu, i, s$. Note that to obtain the asymptotic formulas of high accuracy by the classical methods it is required that $P_2, P_3, ..., P_n$ be differentiable (see [1]). To obtain the asymptotic formulas for $L$ we take the operator $L(C)$, where $L(P_2, P_3, ..., P_n)$ is denoted by $L(C)$ when $P_2(x) = C$ and $P_v(x) = 0$ for $v = 3, 4, ..., n$, for an unperturbed operator and $L - L(C)$ for a perturbation. One can easily verify that the eigenvalues and the normalized eigenfunctions of $L(C)$ are

$$\mu_{k,j} = (2\pi k i)^n + \mu_j (2\pi k i)^{n-2}, \ \Phi_{k,j}(x) = v_j e^{i 2\pi k x} \quad (12)$$

for $k \in \mathbb{Z}$, $j = 1, 2, ..., m$. Since boundary condition (2) is self-adjoint, we have
$(L(C))^* = L(C^*)$. Therefore the eigenfunction $\widetilde{\Phi}_{k,j}(x)$, where $j = 1, 2, ..., m$ and $k \in \mathbb{Z}$, of $(L(C))^*$ corresponding to the eigenvalue $\overline{\mu_{k,j}}$ and satisfying $(\Phi_{k,j}, \widetilde{\Phi}_{k,j}) = 1$ is

$$\widetilde{\Phi}_{k,j}(x) = w_j e^{i2\pi kx}. \tag{13}$$

To prove the asymptotic formulas for the eigenvalue $\lambda_{k,j}$ and for the corresponding normalized eigenfunction $\Psi_{k,j}(x)$ of $L$ we use the formula

$$(\lambda_{k,j} - \mu_{p,s})\left(\Psi_{k,j}, \widetilde{\Phi}_{p,s}\right) = \left((P_2 - C)\Psi_{k,j}^{(n-2)}, \widetilde{\Phi}_{p,s}\right) + \sum_{\nu=3}^{n}\left(P_\nu \Psi_{k,j}^{(n-\nu)}, \widetilde{\Phi}_{p,s}\right) \tag{14}$$

which can be obtained from $L\Psi_{k,j}(x) = \lambda_{k,j}\Psi_{k,j}(x)$ by multiplying scalarly by $\widetilde{\Phi}_{p,s}(x)$. Moreover, we use the following obvious proposition about the system of the eigenfunctions of the operator $L(C)$. We do not consider the statements of the proposition as new. However we could not find a proper reference to all assertions of the proposition and decided to present a short proof here.

**Proposition 1** *If the eigenvalues $\mu_1, \mu_2, ..., \mu_m$ of the matrix $C$ are simple, then for $f \in L_2^m[0, 1]$ the followings hold*

$$f(x) = \sum_{p \in \mathbb{Z};\ q=1,2,...,m} \left(f, \widetilde{\Phi}_{p,q}\right) \Phi_{p,q}(x), \tag{15}$$

$$\| V \|^{-2} \| f \|^2 \leq \sum_{p \in \mathbb{Z};\ q=1,2,...,m} \left|\left(f, \widetilde{\Phi}_{p,q}\right)\right|^2 \leq \| W \|^2 \| f \|^2, \tag{16}$$

*where $V$ and $W$ are the matrices with the columns $v_1, v_2, ..., v_m$ and $w_1, w_2, ..., w_m$.*

**Proof.** Since $v_1, v_2, ..., v_m$ is a basis of $\mathbb{C}^m$ and $\{e^{i2\pi kx} : k \in \mathbb{Z}\}$ is an orthonormal basis of $L_2[0, 1]$, the system

$$\{\Phi_{p,q}, : p \in \mathbb{Z}, q = 1, 2, ..., m\} \tag{17}$$

is a basis of $L_2^m[0, 1]$. Moreover, the sequence $\{\widetilde{\Phi}_{p,q}, : p \in \mathbb{Z}, q = 1, 2, ..., m\}$ is biorthogonal to (17). Therefore, we have (15).

Using the obvious equalities $\Phi_{p,q}(x) = Ve_q e^{i2\pi px}$, $\widetilde{\Phi}_{p,q}(x) = We_q e^{i2\pi px}$, $\| V^* \| = \| V \|$, $\| W^* \| = \| W \|$ and taking into account that the sequence

$$\{e_q e^{i2\pi px} : p \in \mathbb{Z}, q = 1, 2, ..., m\}$$

is an orthonormal basis of $L_2^m[0, 1]$, one can readily see that

$$\sum_{\substack{p \in \mathbb{Z};\\ q=1,2,...,m}} |(f, \Phi_{p,q})|^2 = \sum_{\substack{p \in \mathbb{Z};\\ q=1,2,...,m}} |\left(V^*f, e_q e^{i2\pi px}\right)|^2 = \| V^*f \|^2 \leq \| V \|^2 \| f \|^2, \tag{18}$$

$$\sum_{\substack{p \in \mathbb{Z};\\ q=1,2,...,m}} \left|\left(f, \widetilde{\Phi}_{p,q}\right)\right|^2 = \sum_{\substack{p \in \mathbb{Z};\\ q=1,2,...,m}} |\left(W^*f, e_q e^{i2\pi px}\right)|^2 = \| W^*f \|^2 \leq \| W \|^2 \| f \|^2. \tag{19}$$

On the other hand, it follows from (15) and from the equality

$$f(x) = \sum_{p \in \mathbb{Z};\ q=1,2,...,m} (f, \Phi_{p,q}) \widetilde{\Phi}_{p,q}(x)$$

that
$$\| f \|^2 \leq \sum_{p \in \mathbb{Z};\ q=1,2,...,m} \left|\left(f, \widetilde{\Phi}_{p,q}\right)\right| | (f, \Phi_{p,q}) |.$$

Now using the Schwarz inequality and (18), we get

$$\| f \|^2 \leq \left( \sum_{p \in \mathbb{Z};\ q=1,2,...,m} \left|\left(f, \widetilde{\Phi}_{p,q}\right)\right|^2 \right)^{\frac{1}{2}} \| V \| \| f \|. \tag{20}$$

Inequalities (20) and (19) imply (16) ∎

Formula (11) shows that if $| k | \gg 1$, then the eigenvalue $\lambda_{k,j}$ of $L$ lies far from the eigenvalues $\mu_{p,s}$ for $p \neq \pm k$, namely

$$|\lambda_{k,j} - \mu_{p,s}| > (||k| - |p||)(|k| + |p|)^{n-1}.$$

Using this one can easily verify that

$$\sum_{p:p>d} \frac{|p|^{n-\nu}}{|\lambda_{k,j} - \mu_{p,s}|} = O\left(\frac{1}{d^{\nu-1}}\right), \quad \forall d \geq 2 \mid k \mid, \tag{21}$$

$$\sum_{p:p \neq \pm k} \frac{|p|^{n-\nu}}{|\lambda_{k,j} - \mu_{p,s}|} = O\left(\frac{\ln |k|}{k^{\nu-1}}\right), \tag{22}$$

$$\sum_{p:p \neq \pm k} \frac{|k|^{2n-4}}{|\lambda_{k,j} - \mu_{p,s}|^2} = O\left(\frac{1}{k^2}\right), \tag{23}$$

where $| k | \gg 1$, $\nu \geq 2$.

To estimate the right-hand side of (14) we use (21)-(23) and the following lemma from [14] ( see Lemma 1 of [14]).

*Lemma. Let $\Psi_{k,j,t}(x)$ be the normalized eigenfunction of the operator $L_t$, generated by (1) and the t-periodic boundary conditions*

$$y^{(\nu)}(1) = e^{it} y^{(\nu)}(0), \quad \nu = 0, 1, ..., (n-1)$$

*corresponding to the eigenvalue $\lambda_{k,j}(t) = (2\pi ki + it)^n + O\left(k^{n-1-\frac{1}{2m}}\right)$. Then*

$$\sup_{x \in [0,1]} \left|\Psi_{k,j,t}^{(\nu)}(x)\right| = O(k^\nu) \tag{24}$$

*for $\nu = 0, 1, ..., n-2$. Equality (24) is uniform with respect to t in $[-\frac{\pi}{2}, \frac{3\pi}{2})$.*

It follows from this lemma that

$$\sup_{x \in [0,1]} \left|\Psi_{k,j}^{(\nu)}(x)\right| = O(k^\nu) \tag{25}$$

for $\nu = 0, 1, ..., n-2$ and for $j = 1, 2, ..., m$. Therefore

$$\left((P_2 - C)\Psi_{k,j}^{(n-2)}, \widetilde{\Phi}_{p,s}\right) = O(k^{n-2}), \tag{26}$$

$$\left(P_\nu \Psi_{k,j}^{(n-\nu)}, \widetilde{\Phi}_{p,s}\right) = O(k^{n-\nu}) \tag{27}$$

for all $j, p, s$ and for $\nu = 3, 4, ..., n$. Now (26), (27) and (14) imply that there exist constants $c_1 > 0$ and $N \gg 1$ such that

$$\left|\left(\Psi_{k,j}, \widetilde{\Phi}_{p,q}\right)\right| \leq \frac{c_1 | k |^{n-2}}{| \lambda_{k,j} - \mu_{p,s} |}, \quad \forall p \neq \pm k, \ \forall | k | \geq N, \ \forall j, s. \tag{28}$$

To obtain the asymptotic formulas we use (14), (26)-(28) and Proposition 1.

## 2  Main Results

To prove the main results, first, we prove the following lemma.

**Lemma 1** *The equalities*

$$\left((P_2 - C)\Psi_{k,j}^{(n-2)}, \widetilde{\Phi}_{k,s}\right) \qquad (29)$$
$$= \sum_{q=1,2,\ldots m;\ p=\pm k} (2\pi pi)^{n-2} \left((P_2 - C)\Phi_{p,q}, \widetilde{\Phi}_{k,s}\right) \left(\Psi_{k,j}, \widetilde{\Phi}_{p,q}\right) + O(k^{n-3}\ln|k|),$$

$$\left((P_2 - C)\Phi_{k,q}, \widetilde{\Phi}_{k,s}\right) = 0, \quad \left((P_2 - C)\Phi_{-k,q}(x), \widetilde{\Phi}_{k,s}\right) = b_{s,q,2k} \qquad (30)$$

*hold for all $q$ and $s$.*

**Proof.** Using the integration by parts and (28), we get

$$\left|\left(\Psi_{k,j}^{(n-2)}, \widetilde{\Phi}_{p,q}\right)\right| = |(2\pi p)^{n-2}\left(\Psi_{k,j}, \widetilde{\Phi}_{p,q}\right)| \leq \frac{c_1 |2\pi p|^{n-2} |k|^{n-2}}{|\lambda_{k,j} - \mu_{p,q}|} \qquad (31)$$

for $p \neq \pm k$, $|k| \geq N$. This and (21) imply that there exists a constant $c_2$ such that

$$\sum_{p:|p|>d} |\left(\Psi_{k,j}^{(n-2)}, \widetilde{\Phi}_{p,q}\right)| < \frac{c_2 |k|^{n-2}}{d}$$

for $d \geq 2|k|$. Hence the decomposition of $\Psi_{k,j}^{(n-2)}$ by basis (17) has the form

$$\Psi_{k,j}^{(n-2)}(x) = \sum_{|p|\leq d;\ q=1,2,\ldots,m} (2\pi pi)^{n-2} \left(\Psi_{k,j}, \widetilde{\Phi}_{p,q}\right) \Phi_{p,q}(x) + g_d(x), \qquad (32)$$

where

$$\sup_{x\in[0,1]} |g_d(x)| < \frac{c_2 |k|^{n-2}}{d}.$$

Using (32) in the left-hand side of (29) and letting $d$ tend to $\infty$, we obtain

$$\left((P_2 - C)\Psi_{k,j}^{(n-2)}, \widetilde{\Phi}_{k,s}\right) = \sum_{q=1,2,\ldots m;\ p\in\mathbb{Z}} (2\pi pi)^{n-2} \left((P_2 - C)\Phi_{p,q}, \widetilde{\Phi}_{k,s}\right) \left(\Psi_{k,j}, \widetilde{\Phi}_{p,q}\right). \qquad (33)$$

Since

$$(2\pi pi)^{n-2} \left((P_2 - C)\Phi_{p,q}(x), \widetilde{\Phi}_{k,s}\right) = O(p^{n-2}),$$

it follows from (28) and (22) that

$$\sum_{q=1,2,\ldots,m;\ p\in\mathbb{Z}\setminus\{k,-k\}} (2\pi pi)^{n-2} \left((P_2 - C)\Phi_{p,q}(x), \widetilde{\Phi}_{k,s}\right) \left(\Psi_{k,j}, \widetilde{\Phi}_{p,q}\right) = O(k^{n-3}\ln|k|).$$

This and (33) imply (29).

Using (12) and (13), we obtain

$$\left((P_2 - C)\Phi_{k,q}, \widetilde{\Phi}_{k,s}\right) = \int_0^1 \langle(P_2(x) - C)v_q, w_s\rangle dx, \qquad (34)$$

$$\left((P_2 - C)\Phi_{-k,q}, \widetilde{\Phi}_{k,s}\right) = \int_0^1 \langle(P_2(x) - C)v_q, w_s\rangle e^{-4\pi ikx} dx. \qquad (35)$$

On the other hand, from (6) we have

$$\int_0^1 (P_2(x) - C)dx = 0. \tag{36}$$

Equalities (34) and (36) imply the first equality in (30).

Since $\langle Cv_q, w_s \rangle$ is a constant, we have

$$\int_0^1 \langle Cv_q, w_s \rangle e^{-4\pi ikx} dx = 0.$$

Therefore, the second equality in (30) follows from (35) and (7). ∎

From (29) and (30) we obtain

$$\left((P_2 - C)\Psi_{k,j}^{(n-2)}, \widetilde{\Phi}_{k,s}\right) = (2\pi ki)^{n-2} \sum_{q=1,2,\ldots m} b_{s,q,2k} \left(\Psi_{k,j}, \widetilde{\Phi}_{-k,q}\right) + O(k^{n-3} \ln |k|). \tag{37}$$

This with (27) shows that formula (14) for $p = k$ can be written in the form

$$(\lambda_{k,j} - \mu_{k,s}) \left(\Psi_{k,j}, \widetilde{\Phi}_{k,s}\right) = (2\pi ki)^{n-2} \sum_{q=1,2,\ldots m} b_{s,q,2k} \left(\Psi_{k,j}, \widetilde{\Phi}_{-k,q}\right) + O(k^{n-3} \ln |k|). \tag{38}$$

In the left-hand side of (38) replacing $\widetilde{\Phi}_{k,s}$ by $\widetilde{\Phi}_{-k,s}$ and hence in the right-hand side of (38) replacing $\widetilde{\Phi}_{-k,q}$ by $\widetilde{\Phi}_{k,q}$, we get

$$(\lambda_{k,j} - \mu_{k,s}) \left(\Psi_{k,j}, \widetilde{\Phi}_{-k,s}\right) = (2\pi ki)^{n-2} \sum_{q=1,2,\ldots m} b_{s,q,-2k} \left(\Psi_{k,j}, \widetilde{\Phi}_{k,q}\right) + O(k^{n-3} \ln |k|). \tag{39}$$

Using (38), (39) and (7) one can readily see that there exists a constant $c_3$ such that

$$\left| (\lambda_{k,j} - \mu_{k,s}) \left(\Psi_{k,j}, \widetilde{\Phi}_{\pm k,s}\right) \right| < c_3 k^{n-2} (b_{\pm 2k} + |k|^{-1} \ln |k|). \tag{40}$$

Let $\varepsilon_k = 2mc_3 k^{n-2}(b_{2k} + b_{-2k} + |k|^{-1} \ln |k|) \| V \|$.

**Theorem 1** *Suppose the eigenvalues $\mu_1, \mu_2, \ldots, \mu_m$ of the matrix $C$ are simple. Then:*

*(a) There exist a number $N_0 \geq N$, where $N$ is defined in (28), such that the eigenvalues $\lambda_{k,1}, \lambda_{k,2}, \ldots, \lambda_{k,m}$ of $L$ for $|k| \geq N_0$ lie in the union of the pairwise disjoint disks*

$$U(\mu_{k,1}, \varepsilon_k), \ U(\mu_{k,2}, \varepsilon_k), \ldots, U(\mu_{k,m}, \varepsilon_k), \tag{41}$$

*where $U(\mu, c) = \{\lambda \in \mathbb{C} : |\lambda - \mu| < c\}$.*

*(b) For each $j$ and for $|k| \geq N_0$ the disk $U(\mu_{k,j}, \varepsilon_k)$ contains precisely 2 eigenvalues (counting multiplicity), denoted by $\lambda_{k,j}$ and $\lambda_{-k,j}$. If $q \neq j$, then the equality*

$$\left(\Psi_{k,j}, \widetilde{\Phi}_{\pm k,q}\right) = O(b_{\pm 2k}) + O(|k|^{-1} \ln |k|) \tag{42}$$

*holds for any eigenfunction $\Psi_{k,j}$ corresponding to any of the eigenvalues $\lambda_{k,j}$ and $\lambda_{-k,j}$.*

**Proof.** (a) Suppose to the contrary that $\lambda_{k,j} \notin U(\mu_{k,s}, \varepsilon_k)$ for all $s$. Then we have

$$|\lambda_{k,j} - \mu_{k,s}| \geq \varepsilon_k, \ \forall s$$

This and (40) imply that

$$\left| \left(\Psi_{k,j}, \widetilde{\Phi}_{\pm k,s}\right) \right| < \frac{1}{2m \| V \|}.$$

Then
$$\sum_{p=\pm k,\ s=1,2,\dots,m} \left|\left(\Psi_{k,j}, \widetilde{\Phi}_{p,s}\right)\right|^2 < \frac{1}{2} \| V \|^{-2}.$$

On the other hand, it follows from (28) and (23) that
$$\sum_{p\neq \pm k,\ s=1,2,\dots,m} \left|\left(\Psi_{k,j}, \widetilde{\Phi}_{p,s}\right)\right|^2 = O(k^{-2}). \tag{43}$$

The last 2 relations and the equality $\| \Psi_{k,j} \|= 1$ contradict (16).

It follows from the definitions of $b_k$, $\varepsilon_k$ and $\mu_{k,j}$ that
$$\lim_{k\to\infty} b_{\pm 2k} = 0,\ \varepsilon_k = o(k^{n-2}),\ | \mu_{k,j} - \mu_{k,q} |\geq a(2\pi k)^{n-2} \tag{44}$$

for all $q \neq j$, where
$$a = \min_{k\neq s} | \mu_k - \mu_s |.$$

Therefore the disks in (41) are pairwise disjoint.

(b) Consider the following family of operators
$$L_\varepsilon = L(C) + \varepsilon(L - L(C)),\ 0 \leq \varepsilon \leq 1.$$

The formula (14) for the operator $L_\varepsilon$ has the form
$$(\lambda_{k,j,\varepsilon} - \mu_{p,s})\left(\Psi_{k,j,\varepsilon}, \widetilde{\Phi}_{p,s}\right) = \varepsilon\left((P_2 - C)\Psi_{k,j,\varepsilon}^{(n-2)}, \widetilde{\Phi}_{p,s}\right) + \varepsilon\sum_{\nu=3}^{n}\left(P_\nu \Psi_{k,j,\varepsilon}^{(n-\nu)}, \widetilde{\Phi}_{p,s}\right),$$

where $\lambda_{k,j,\varepsilon}$ and $\Psi_{k,j,\varepsilon}$ are the eigenvalue and eigenfunction of $L_\varepsilon$. Therefore using this formula instead of (14) and repeating the arguments by which we obtained the proof of the case (a) of Theorem 1, one can see that the assertions of the case (a) of Theorem 1 hold for $L_\varepsilon$. It means that the eigenvalues $\lambda_{k,j,\varepsilon}$ of $L_\varepsilon$ for $| k |\geq N_0$ lie in the union of the disks in (41). Hence the boundary $\partial(U(\mu_{k,j}, \varepsilon_k))$ of the disk $U(\mu_{k,j}, \varepsilon_k)$ lies in the resolvent set of $L_\varepsilon$ for $\varepsilon \in [0,1]$. Therefore, taking into account that the family $L_\varepsilon$ is halomorphic (in the sense of [15]) with respect to $\varepsilon$, we obtain that the number of the eigenvalues of $L_\varepsilon$ lying inside of $\partial(U(\mu_{k,j}, \varepsilon_k))$ are the same for all $\varepsilon \in [0,1]$. Since $L_0 = L(C)$ and $L(C)$ has only one eigenvalue $\mu_{k,j}$ of multiplicity 2 in the disks $U(\mu_{k,j}, \varepsilon_k)$, the operator $L$ has two eigenvalues (counting multiplicity) in the disk $U(\mu_{k,1}, \varepsilon_k)$. Using (44) and the inclusion $\lambda_{k,j} \in U(\mu_{k,j}, \varepsilon_k)$ we see that
$$| \lambda_{k,j} - \mu_{k,q} |> \frac{1}{2}a(2\pi k)^{n-2},\ \forall q \neq j.$$

Therefore (42) follows from (40). ∎

Using (42) in (38) and (39) and then taking into account (7), we obtain
$$(\lambda_{k,j} - \mu_{k,j})\left(\Psi_{k,j}, \widetilde{\Phi}_{k,j}\right) = (2\pi k i)^{n-2}\left(b_{j,j,2k}\left(\Psi_{k,j}, \widetilde{\Phi}_{-k,j}\right) + O(b_{2k}b_{-2k}) + O\left(\frac{\ln|k|}{k}\right)\right),$$

$$(\lambda_{k,j} - \mu_{k,j})\left(\Psi_{k,j}, \widetilde{\Phi}_{-k,j}\right) = (2\pi k i)^{n-2}\left(b_{j,j,-2k}\left(\Psi_{k,j}, \widetilde{\Phi}_{k,j}\right) + O(b_{2k}b_{-2k}) + O\left(\frac{\ln|k|}{k}\right)\right).$$

Dividing both sides of these equalities by $(2i\pi k)^{n-2}$, we get
$$\left(\Lambda_{k,j} - (2i\pi k)^2 - \mu_j\right) u_{k,j} = b_{j,j,2k} v_{k,j} + O(d_k), \tag{45}$$

$$(\Lambda_{k,j} - (2i\pi k)^2 - \mu_j)v_{k,j} = b_{j,j,-2k}u_{k,j} + O(d_k), \tag{46}$$

where

$$\Lambda_{k,j} = \frac{\lambda_{k,j}}{(2i\pi k)^{n-2}}, \ u_{k,j} = \left(\Psi_{k,j}, \widetilde{\Phi}_{k,j}\right), \ v_{k,j} = \left(\Psi_{k,j}, \widetilde{\Phi}_{-k,j}\right), \tag{47}$$

$$d_k = \max\{b_{2k}b_{-2k}, |k|^{-1}\ln|k|\}, \ \lambda_{k,j} \in U(\mu_{k,j}, \varepsilon_k). \tag{48}$$

Using (43), (42) and Proposition 1 and taking into account that
$(\Phi_{k,j}, \Phi_{-k,j}) = 0$, $\|\Psi_{k,j}\| = 1$, $\|\Phi_{\pm k,j}\| = 1$, we obtain

$$\Psi_{k,j} = u_{k,j}\Phi_{k,j} + v_{k,j}\Phi_{-k,j} + O(b_{2k}) + O(b_{-2k}) + O(|k|^{-1}\ln|k|), \ |u_{k,j}|^2 + |v_{k,j}|^2 = 1 + o(1). \tag{49}$$

Now, using (45)-(49), we obtain asymptotic formulas.

**Theorem 2** *Suppose the eigenvalues $\mu_1, \mu_2, ..., \mu_m$ of the matrix $C$ are simple. Let $\lambda_{k,j}$ be an eigenvalues of $L$ lying in $U(\mu_{k,j}, \varepsilon_k)$. If condition (8) holds, then there exist numbers $c_4 > 0$ and $N_1 \geq N_0$, where $N_0$ is defined in Theorem 1, such that:*

*(a) The eigenvalue $\lambda_{k,j}$ for $|k| \geq N_1$ lies in $U_{-k,j} \cup U_{k,j}$, where*

$$U_{\pm k,j} = \{\lambda \in \mathbb{C}: \ |\lambda - h_{\pm k,j}| < c_4 k^{n-2}\gamma_{2k}d_k\}, \quad U_{k,j} \cap U_{-k,j} = \emptyset, \tag{50}$$

$$h_{\pm k,j} = (i2\pi k)^n + \mu_j(2\pi ki)^{n-2} \pm (2\pi ki)^{n-2}q_{2k}, \ q_k = (b_{j,j,k}b_{j,j,-k})^{\frac{1}{2}},$$

$$\gamma_k = \max\left\{\left(\frac{|b_{j,j,k}|}{|b_{j,j,-k}|}\right)^{\frac{1}{2}}, \left(\frac{|b_{j,j,-k}|}{|b_{j,j,k}|}\right)^{\frac{1}{2}}\right\}. \tag{51}$$

*(b) The geometrical multiplicity of the eigenvalue $\lambda_{k,j}$ for $|k| \geq N_1$ is 1. If $\lambda_{k,j}$ lies in $U_{\pm k,j}$ then any eigenfunction $\Psi_{k,j}$ of $L$ corresponding to $\lambda_{k,j}$ satisfies*

$$\Psi_{k,j} = \left(1 + |\alpha_{2k,j}|^2\right)^{-\frac{1}{2}}(\Phi_{k,j} + \alpha_{\pm 2k,j}\Phi_{-k,j}) + O(b_{2k}) + O(b_{-2k}) + O(|k|^{-1}\ln|k|), \tag{52}$$

*where*

$$\alpha_{\pm k,j} = \frac{\pm q_k}{b_{j,j,k}}(1 + o(1)). \tag{53}$$

**Proof.** (*a*) We use the following equalities that easily follow from (8) and from the definitions of $d_k, q_k, \gamma_k$, (see (48), (50), (51))

$$\frac{d_k}{b_{j,j,\pm 2k}} = o(1), \ \frac{\gamma_{2k}\ln|k|}{kq_{2k}} = o(1), \ \frac{\gamma_{2k}d_k}{q_{2k}} = o(1). \tag{54}$$

The last equality in (54) implies the second relation in (50).

Since (49) holds, at least one of the numbers $|u_{k,j}|, |v_{k,j}|$ is greater than $\frac{1}{2}$ and the inequalities $|u_{k,j}| < 2, |v_{k,j}| < 2$ are satisfied. Therefore, at least one of the following relations holds

$$|u_{k,j}| \sim 1, \quad |v_{k,j}| \sim 1. \tag{55}$$

Assume that the first relation of (55) holds. Then dividing (45) by $u_{k,j}$, we get

$$(\Lambda_{k,j} - (2i\pi k)^2 - \mu_j) = b_{j,j,2k}\frac{v_{k,j}}{u_{k,j}} + O(d_k). \tag{56}$$

Now we estimate $\frac{v_{k,j}}{u_{k,j}}$ as follows: multiply (45) and (46) by $v_{k,j}$ and $u_{k,j}$ respectively and take the difference to get

$$b_{j,j,2k}v_{k,j}^2 = b_{j,j,-2k}u_{k,j}^2 + O(d_k), \quad \left(\frac{v_{k,j}}{u_{k,j}}\right)^2 = \frac{b_{j,j,-2k}}{b_{j,j,2k}}\left(1 + O\left(\frac{d_k}{b_{j,j,-2k}}\right)\right).$$

Now using the first equality of (54), we obtain

$$\left(\frac{v_{k,j}}{u_{k,j}}\right) = \left(\frac{b_{j,j,-2k}}{b_{j,j,2k}}\right)^{\frac{1}{2}}\left(1 + O\left(\frac{d_k}{b_{j,j,-2k}}\right)\right).$$

This and (51) imply

$$b_{j,j,2k}\frac{v_{k,j}}{u_{k,j}} = q_{2k} + O(\gamma_{2k}d_k).$$

Using this in (56), and taking into account that $\gamma_{2k} \geq 1$ (see (51)), we get

$$\mid \Lambda_{k,j} - (2i\pi k)^2 - \mu_j \mid = \mid q_{2k} \mid + O(\gamma_{2k}d_k). \tag{57}$$

If the second relation of (55) holds, then in the same way we obtain (57). Now the definition of $\Lambda_{k,j}$ (see (47)) and (57) imply the proof of $(a)$.

$(b)$ If $\lambda_{k,j}$ lies in $U_{\pm k,j}$, then by the definitions of $U_{\pm k,j}$ and $\Lambda_{k,j}$, we have

$$\Lambda_{k,j} = (i2\pi k)^2 + \mu_j \pm q_{2k} + O(\gamma_{2k}d_k). \tag{58}$$

Substituting (58) into (45) and (46), we obtain the equalities

$$\pm q_{2k}u_{k,j} = b_{j,j,2k}v_{k,j} + O(\gamma_{2k}d_k), \quad \pm q_{2k}v_{k,j} = b_{j,j,-2k}u_{k,j} + O(\gamma_{2k}d_k).$$

Using the first equality if the first relation of (55) holds and using the second equality if the second relation of (55) holds, and taking into account (54), we see that

$$\frac{v_{k,j}}{u_{k,j}} = \frac{\pm q_{2k}}{b_{j,j,2k}}(1 + o(1)). \tag{59}$$

Now (59), (47) and (49) imply (52) and (53). If there are two linearly independent eigenfunctions corresponding to $\lambda_{k,j}$, then one can find two orthogonal eigenfunctions satisfying (52), which is impossible ∎

Now we prove that the eigenvalues $\lambda_{k,j}$ for large value of $k$ are simple and in each of the disks $U_{-k,j}$ and $U_{k,j}$ defined in (50) there exists unique eigenvalue of $L$. For this we consider the following family of operators

$$B_\varepsilon = S + \varepsilon(L - S), \ 0 \leq \varepsilon \leq 1, \tag{60}$$

where $S$ is the operator generated by (2) and by the differential expression

$$y^{(n)} + (C + (b_{j,j,2k}e^{i4\pi kx} + b_{j,j,-2k}e^{-i4\pi kx})I)y^{(n-2)}, \tag{61}$$

$I$ is $m \times m$ unit matrix. We denote by $\lambda_{k,j,\varepsilon}$ and $\Psi_{k,j,\varepsilon}$ the eigenvalue and eigenfunction of $B_\varepsilon$. Note that this notations were used for the eigenvalue and eigenfunction of $L_\varepsilon$ in the proof of Theorem 1. Here, for simplicity of notation, we use the same symbols.

**Lemma 2** *Suppose the eigenvalues $\mu_1, \mu_2, ..., \mu_m$ of the matrix $C$ are simple. Let $\lambda_{k,j,0}$ be an eigenvalue of $S$ lying in the disk $U_{\pm k,j}$ defined in (50). Then any normalized eigenfunctions $\Psi_{k,j,0}(x)$ and $\widetilde{\Psi}_{k,j,0}(x)$ of $S$ and $S^*$ corresponding to the eigenvalues $\lambda_{k,j,0}$ and $\overline{\lambda_{k,j,0}}$ respectively satisfy*

$$\Psi_{k,j,0}(x) = \left(1 + \mid \alpha_{2k,j} \mid^2\right)^{-\frac{1}{2}}\left(\Phi_{k,j}(x) + \alpha_{\pm 2k,j}\Phi_{-k,j}(x)\right) + O\left(\frac{1}{k}\right), \tag{62}$$

$$\widetilde{\Psi}_{k,j,0}(x) = \left(1+ \mid \widetilde{\alpha}_{2k,j} \mid^2\right)^{-\frac{1}{2}} \left(\Phi_{k,j}(x) + \widetilde{\alpha}_{\pm 2k,j}\Phi_{-k,j}(x)\right) + O\left(\frac{1}{k}\right), \qquad (63)$$

where $\alpha_{k,j}$ is defined in (53) and $\widetilde{\alpha}_{\pm k,j} = \frac{\pm \bar{q}_k}{b_{j,j,-k}}(1+o(k))$.

**Proof.** The proof of (62) is similar to the proof of (52). Formula (52) for the eigenfunction $\Psi_{k,j}$ of $L$ is obtained from the formulas (45), (46) and (49) for $L$. By the same argument we can establish (62) provided suitable formulas like (45), (46) and (49) are obtained for $S$. Formula (14) for $S$ has the form

$$(\lambda_{k,j,0} - \mu_{p,s})\left(\Psi_{k,j,0}, \widetilde{\Phi}_{p,s}\right) = b_{j,j,2k}\left(\Psi_{k,j,0}^{(n-2)}, \widetilde{\Phi}_{p-2k,s}\right) + b_{j,j,-2k}\left(\Psi_{k,j,0}^{(n-2)}, \widetilde{\Phi}_{p+2k,s}\right), \quad (64)$$

since $S$ can be obtained from $L$ by taking

$$P_2(x) = (C + (b_{j,j,2k}e^{i4\pi kx} + b_{j,j,-2k}e^{-i4\pi kx})I), \ P_v(x) = 0, \ \forall v > 2.$$

In (64) replacing $p$ by $3k$, then dividing by $\lambda_{k,j} - \mu_{3k,s}$ and then using the obvious relations

$$\left(\Psi_{k,j,0}^{(n-2)}, \widetilde{\Phi}_{q,s}\right) = (2\pi qi)^{n-2}\left(\Psi_{k,j,0}, \widetilde{\Phi}_{q,s}\right), \ \mid \lambda_{k,j} - \mu_{3k,s} \mid > k^n \qquad (65)$$

for $\mid k \mid \gg 1$, we obtain

$$\left(\Psi_{k,j,0}, \widetilde{\Phi}_{\pm 3k,s}\right) = O(k^{-2}). \qquad (66)$$

Now in (64) replace $p$ and $s$ by $\pm k$ and $j$ respectively, use (65), (66) and the notations

$$\Lambda_{k,j,0} = \frac{\lambda_{k,j,0}}{(2i\pi k)^{n-2}}, \ u_{k,j,0} = \left(\Psi_{k,j,0}, \widetilde{\Phi}_{k,j}\right), \ v_{k,j,0} = \left(\Psi_{k,j,0}, \widetilde{\Phi}_{-k,j}\right), \qquad (67)$$

to get the equalities

$$(\Lambda_{k,j,0} - (2i\pi k)^2 - \mu_j)u_{k,j,0} = b_{j,j,2k}v_{k,j,0} + O(k^{-2}), \qquad (68)$$

$$(\Lambda_{k,j} - (2i\pi k)^2 - \mu_j)v_{k,j,0} = b_{j,j,-2k}u_{k,j,0} + O(k^{-2}). \qquad (69)$$

Equalities (68) and (69) are the analog of (45) and (46) for $S$.

Now, we obtain the analog of (49) as follows. In (64) replace $p$ by $k$ and then by $-k$, use (65) and (66), to get the equalities

$$(\lambda_{k,j,0} - \mu_{k,s})\left(\Psi_{k,j,0}, \widetilde{\Phi}_{k,s}\right) = (2\pi ki)^{n-2}b_{j,j,2k}\left(\Psi_{k,j,0}, \widetilde{\Phi}_{-k,s}\right) + O(k^{n-4}). \qquad (70)$$

$$(\lambda_{k,j,0} - \mu_{k,s})\left(\Psi_{k,j,0}, \widetilde{\Phi}_{-k,s}\right) = (2\pi ki)^{n-2}b_{j,j,-2k}\left(\Psi_{k,j,0}, \widetilde{\Phi}_{k,s}\right) + O(k^{n-4}). \qquad (71)$$

Since $\lambda_{k,j,0} \in U_{\pm k,j}$, it follows from the inequality in (44) that

$$\mid \lambda_{k,j,0} - \mu_{k,s} \mid > \frac{1}{2}a(2\pi k)^{n-2}, \ \forall s \neq j. \qquad (72)$$

Formulas (70)-(72) with the first equality of (44) yield

$$\left(\Psi_{k,j,0}, \widetilde{\Phi}_{\pm k,s}\right) = O(k^{-2}), \ \forall s \neq j. \qquad (73)$$

This formula and the formula which is obtained from (43) by replacing $\Psi_{k,j}$ with $\Psi_{k,j,0}$ imply that the expansion of $\Psi_{k,j,0}$ has the form

$$\Psi_{k,j,0}(x) = u_{k,j,0}\Phi_{k,j}(x) + v_{k,j,0}\Phi_{-k,j}(x) + h(x), \qquad (74)$$

where
$$\|h(x)\| = O(k^{-1}), \quad |u_{k,j,0}|^2 + |v_{k,j,0}|^2 = 1 + O(k^{-2}). \tag{75}$$

Instead of (45), (46), (49) using (68), (69), (74), (75) and arguing as in the proof of (52), we obtain (62).

To prove (63), we use the formula
$$\left(\overline{\lambda_{k,j,0}} - \overline{\mu_{p,s}}\right)\left(\widetilde{\Psi}_{k,j,0}, \Phi_{p,s}\right) = \overline{b_{j,j,2k}}\left(\widetilde{\Psi}_{k,j,0}^{(n-2)}, \Phi_{p+2k,s}\right) + \overline{b_{j,j,-2k}}\left(\widetilde{\Psi}_{k,j,0}^{(n-2)}, \Phi_{p-2k,s}\right) \tag{76}$$

which can be obtained from
$$S^*\widetilde{\Psi}_{k,j,0}(x) = \overline{\lambda_{k,j,0}}\widetilde{\Psi}_{k,j,0}(x) \tag{77}$$

by multiplying by $\Phi_{p,s}(x)$ and using
$$\left(S^*\widetilde{\Psi}_{k,j,0}, \Phi_{p,s}\right) = \left(\widetilde{\Psi}_{k,j,0}, S\Phi_{p,s}\right) = \overline{\mu_{p,s}}\left(\widetilde{\Psi}_{k,j,0}, \Phi_{p,s}\right) +$$
$$\overline{b_{j,j,2k}}\left(\widetilde{\Psi}_{k,j,0}^{(n-2)}, \Phi_{p+2k,s}\right) + \overline{b_{j,j,-2k}}\left(\widetilde{\Psi}_{k,j,0}^{(n-2)}, \Phi_{p-2k,s}\right).$$

Instead of (64) using (76) and arguing as in the proof of (62) we obtain (63) ∎

**Theorem 3** *Let the eigenvalues $\mu_1, \mu_2, ..., \mu_m$ of the matrix $C$ be simple. If (8) holds, then:*
*(a) The eigenvalues $\lambda_{k,j}$ for $|k| \geq N_1$ are simple and consist of 2 sequences $\{\lambda_{k,j} : k \geq N_1\}$ and $\{\lambda_{-k,j} : k \geq N_1\}$ satisfying*
$$\lambda_{\pm k,j} = (i2\pi k)^n + \mu_j(2\pi ki)^{n-2} \pm (2\pi ki)^{n-2}q_{2k} + O\left(k^{n-3}\gamma_{2k}\ln|k|\right),$$

*where $j = 1, 2, ..., m$ and $N_1$ is defined in Theorem 2. The normalized eigenfunction $\Psi_{\pm k,j}(x)$ of $L$ corresponding to the eigenvalue $\lambda_{\pm k,j}$ satisfies*
$$\Psi_{\pm k,j}(x) = (1+|\alpha_{\pm 2k,j}|^2)^{-\frac{1}{2}}(\Phi_{k,j}(x) + \alpha_{\pm 2k,j}\Phi_{-k,j}(x)) + O(b_{2k}^2) + O(b_{-2k}^2) + O\left(\frac{\ln|k|}{k}\right),$$

*where $q_k$, $\gamma_k$ and $\alpha_{\pm k,j}$ are defined in Theorem 2.*
*(b) The root functions of $L$ form a Riesz basis in $L_2^m[0,1]$ if and only if $b_{j,j,2k} \sim b_{j,j,-2k}$ for all $j = 1, 2, ..., m$.*

**Proof.** (a) Formula (14) for the operator $B_\varepsilon$ has the form
$$(\lambda_{k,j,\varepsilon} - \mu_{p,s})\left(\Psi_{k,j,\varepsilon}, \widetilde{\Phi}_{p,s}\right) = \varepsilon\left((P_2 - C)\Psi_{k,j,\varepsilon}^{(n-2)}, \widetilde{\Phi}_{p,s}\right) + \varepsilon\sum_{\nu=3}^{n}\left(P_\nu\Psi_{k,j,\varepsilon}^{(n-\nu)}, \widetilde{\Phi}_{p,s}\right) +$$
$$(1-\varepsilon)\left(b_{j,j,2k}\left(\Psi_{k,j,\varepsilon}^{(n-2)}, \widetilde{\Phi}_{p-2k,s}\right) + b_{j,j,-2k}\left(\Psi_{k,j,\varepsilon}^{(n-2)}, \widetilde{\Phi}_{p+2k,s}\right)\right).$$

Instead of (14) using this formula and repeating the proof of Theorem 2, one can see that the assertions of Theorem 2 hold for the operator $B_\varepsilon$. Thus
$$\{\lambda_{k,j,\varepsilon}, \lambda_{-k,j,\varepsilon}\} \subset U_{k,j} \cup U_{-k,j}, \ U_{k,j} \cap U_{-k,j} = \emptyset, \ \forall k \geq N_1, \ \forall \varepsilon \in [0,1]. \tag{78}$$

Now using Lemma 2, we prove that the eigenvalue $\lambda_{k,j,0}$ of $S$ lying in the disk $U_{\pm k,j}$ for large value of $k$ is simple. Since the geometrical multiplicity of this eigenvalue is 1 (see Theorem 2), we need to prove that there is not associated function corresponding to $\Psi_{k,j,0}(x)$. Suppose to the contrary that there exists an associated function of $S$ corresponding to $\Psi_{k,j,0}(x)$. Then $(\Psi_{k,j,0}, \widetilde{\Psi}_{k,j,0}) = 0$. Therefore, using (62), (63), the definition of $q_k$ (see (50)) and the equalities $(\Phi_{k,j}, \Phi_{-k,j}) = 0$, $\| \Psi_{k,j,0} \| = 1$, $\| \widetilde{\Psi}_{k,j,0} \| = 1$, $\| \Phi_{\pm k,j} \| = 1$, we get

$$2\left(1+ \mid \widetilde{\alpha}_{2k,j} \mid^2\right)^{-\frac{1}{2}} \left(1+ \mid \alpha_{2k,j} \mid^2\right)^{-\frac{1}{2}} = O\left(\frac{1}{k}\right). \tag{79}$$

Let us prove that (79) contradicts (8). It follows from (8) that

$$\left|\frac{1}{b_{j,j,2k}}\right| \ll \left|\frac{k}{\ln k}\right|, \quad \left|\frac{1}{b_{j,j,-2k}}\right| \ll \left|\frac{k}{\ln k}\right| \quad \text{for } k \gg 1.$$

This, (8) and the equalities $q_{2k}^2 = b_{j,j,-2k} b_{j,j,2k}$, $\lim_{k \to \infty} b_{j,j,\pm 2k} = 0$ imply

$$1+ \mid \alpha_{2k,j} \mid^2 < \left|\frac{k}{\ln k}\right|, \quad 1+ \mid \widetilde{\alpha}_{2k,j} \mid^2 < \left|\frac{k}{\ln k}\right| \quad \text{for } k \gg 1$$

which contradicts (79). Thus the eigenvalues $\lambda_{k,j,0}$ of $S$ lying in the disk $U_{\pm k,j}$ is simple. It means that the eigenvalues $\lambda_{k,j,0}$ and $\lambda_{-k,j,0}$ are simple.

Now we prove that in each of the intervals $U_{k,j}$ and $U_{-k,j}$ for $k \geq N_1$ there exists unique eigenvalue of $B_0 = S$. Suppose to the contrary that both eigenvalues $\lambda_{k,j,0}$ and $\lambda_{-k,j,0}$ of $S$ lie in the same interval and, without loss of generality, assume that $\{\lambda_{k,j,0}, \lambda_{-k,j,0}\} \subset U_{k,j}$. Then by Lemma 2 both eigenfunctions $\Psi_{k,j,0}(x)$ and $\Psi_{-k,j,0}(x)$ corresponding to $\lambda_{k,j,0}$ and $\lambda_{-k,j,0}$ respectively satisfy the formula obtained from (62) by replacing $\pm$ with $+$. Similarly, both eigenfunctions $\widetilde{\Psi}_{k,j,0}(x)$ and $\widetilde{\Psi}_{-k,j,0}(x)$ of $S^*$ corresponding to the eigenvalues $\overline{\lambda_{k,j,0}}$ and $\overline{\lambda_{-k,j,0}}$ satisfy the formula obtained from (63) by replacing $\pm$ with $+$. Using this in the equality $(\Psi_{k,j,0}, \widetilde{\Psi}_{-k,j,0}) = 0$, we get (79) which, as proved above, contradicts (8). Thus we proved that in each of the disks $U_{k,j}$ and $U_{-k,j}$ there exists unique eigenvalue of $B_0$. By (78) the boundary $\partial(U_{\pm k,j})$ of the disk $U_{\pm k,j}$ lies in the resolvent set of the operators $B_\varepsilon$ for $\varepsilon \in [0,1]$. Therefore taking into account that the family $B_\varepsilon$ is halomorphic with respect to $\varepsilon$, we obtain that the number of the eigenvalues of $B_\varepsilon$ lying inside of $\partial(U_{\pm k,j})$ are the same for all $\varepsilon \in [0,1]$. Therefore in each of the disks $U_{k,j}$ and $U_{-k,j}$ for $k \geq N_1$ there exists unique eigenvalue of $L$. Thus Theorem 2 implies the case $(a)$ of Theorem 3.

$(b)$ Using the asymptotic formulas for $\Psi_{k,j}$ and $\Psi_{-k,j}$ obtained in the case $(a)$ of Theorem 3, we obtain

$$(\Psi_{k,j}, \Psi_{-k,j}) = \left(1 - \left|\frac{b_{j,j,-2k}}{b_{j,j,2k}}\right|\right)\left(1 + \left|\frac{b_{j,j,-2k}}{b_{j,j,2k}}\right|\right)^{-1} + o(1).$$

This implies that $b_{j,j,2k} \sim b_{j,j,-2k}$ if and only if the following holds:

$$\exists a \in (0,1) \text{ such that } \sup_{k \geq N_1} \mid (\Psi_{k,j}, \Psi_{-k,j}) \mid < a, \quad \forall j = 1, 2, ..., m. \tag{80}$$

It remains to prove that the root functions of $L$ form a Riesz basis if and only if (80) holds. If (80) does not hold, then there exist sequences

$\{k_s : s = 1, 2, ...\}$, $\{a_s \in \mathbb{C} : s = 1, 2, ...\}$, $\{b_s \in \mathbb{C} : s = 1, 2, ...\}$ such that

$$\lim_{s \to \infty} \mid (\Psi_{k_s,j}, \Psi_{-k_s,j}) \mid = 1, \quad \lim_{s \to \infty} \frac{\mid a_s \Psi_{k_s,j} + b_s \Psi_{-k_s,j} \mid^2}{\mid a_s \mid^2 + \mid b_s \mid^2} = 0.$$

This implies that inequality (2.4) in Chapter 6 of [7] (see Theorem 2.1 (N. K. Bari) in Chapter 6 of [7] ) does not hold. Thus, by the Bari Theorem, the root functions of $P$ do not form a Riesz basis.

Now suppose that (80) holds. By Theorem 1 apart from the eigenvalues $\lambda_{k,j}$, where $|k| \geq N_1$, $j = 1, 2, ..., m$, there exist finite eigenvalues $\lambda_1, \lambda_2, ..., \lambda_s$, of the operator $L$. Let $H_k$ be the eigenspace corresponding to the eigenvalue $\lambda_k$ and let $G_k$ be $2m$ dimensional space generated by the eigenfunctions $\Psi_{k,1}, \Psi_{-k,1}, \Psi_{k,2}, \Psi_{-k,2}, ..., \Psi_{k,m}, \Psi_{-k,m}$, where $k \geq N_1$. It is known [16] that the sequence

$$\{H_1, H_2, ..., H_s, G_{N_1}, G_{N_1+1}, ..., \} \tag{81}$$

forms a Riesz basis of subspaces. Let $\varphi_{k,1}, \varphi_{k,2}, ..., \varphi_{k,j_k}$ be an orthonormal basis of the subspace $H_k$. Now we prove that the system

$$(\cup_{k=1}^{s}\{\varphi_{k,1}, \varphi_{k,2}, ..., \varphi_{k,j_k}\}) \cup (\cup_{k \geq N_1}\{\Psi_{k,1}, \Psi_{-k,1}, \Psi_{k,2}, \Psi_{-k,2}, ..., \Psi_{k,m}, \Psi_{-k,m}\})$$

forms an ordinary Riesz basis in $L_2^m[0,1]$. For this we consider the following

$$\Psi = \sum_{k=1}^{s}(a_{k,1}\varphi_{k,1} + a_{k,2}\varphi_{k,2} + ... + a_{k,j_k}\varphi_{k,j_k}) + \sum_{k=N_1}^{N_2}\left(\sum_{j=1,2,...,m}(b_{k,j}\Psi_{k,j} + b_{-k,j}\Psi_{-k,j})\right),$$

where $a_{k,1}, a_{k,2}, ..., a_{k,j_k}$ and $b_{k,j}, b_{-k,j}$ are the complex numbers and $N_2 > N_1$. It follows from (5.24) of section 6 of [7] that

$$\frac{\|\Psi\|^2}{c} \leq \sum_{k=1}^{s}\sum_{j=1}^{j_k}|a_{k,j}|^2 + \sum_{k=N_1}^{N_2}\left\|\sum_{j=1,2,...,m}(b_{k,j}\Psi_{k,j} + b_{-k,j}\Psi_{-k,j})\right\|^2 \leq c\|\Psi\|^2, \tag{82}$$

where $c = \|B\|^2\|B^{-1}\|^2$ and $B$ is a bounded linear invertible operator which transform some orthogonal basis of the subspaces of the space $L_2^m[0,1]$ into basis (81). Using (80) and the asymptotic formulas for $\Psi_{k,j}$ obtained in Theorem 3, taking into account that the normalized eigenvectors $v_1, v_2, ..., v_m$ of the matrix $C$ form a basis of $\mathbb{C}^m$ and all norms are equivalent in the finite dimensional spaces, one can readily see that there exist constants $c_5$ and $c_6$ such that

$$c_5 \sum_{j=1,2,...,m}(|b_{k,j}|^2 + |b_{-k,j}|^2) \leq \left\|\sum_{j=1,2,...,m}(b_{k,j}\Psi_{k,j} + b_{-k,j}\Psi_{-k,j})\right\|^2$$
$$\leq c_6 \sum_{j=1,2,...,m}(|b_{k,j}|^2 + |b_{-k,j}|^2)$$

for $k \geq N_1$. This with (82) implies that inequality (2.4) of the Bari Theorem in Chapter 6 of [7] holds, i.e., the root functions of $L$ form a Riesz basis ∎

Now let us consider the operator $A(P_2, P_3, ..., P_n)$ generated by (1) and the antiperiodic boundary condition (9). Due to the classical investigations [1, chapter 3, theorem 2] all large eigenvalues of $A$ consist of the sequences

$$\{\rho_{k,1} :| k |\geq N\}, \{\rho_{k,2} :| k |\geq N\}, ..., \{\rho_{k,m} :| k |\geq N\},$$

where $N \gg 1$, $k \in \mathbb{Z}$, satisfying the following asymptotic formulas

$$\rho_{k,j} = ((2k+1)\pi i)^n + O\left(k^{n-1-\frac{1}{2m}}\right)$$

for $j = 1, 2, ..., m$. Let $X_{k,j}$ be the eigenfunction of $A$ corresponding to $\rho_{k,j}$. The operator $A(P_2, P_3, ..., P_n)$ is denoted by $A(C)$ when $P_2(x) = C$, $P_v(x) = 0$ for $v = 3, 4, ..., n$. Let $E_{k,s}$ and $\widetilde{E}_{k,s}$ be the eigenfunctions of $A(C)$ and $(A(C))^*$ corresponding to the eigenvalues

$$((2k+1)\pi i)^n + \mu_j((2k+1)\pi i)^{n-2} \text{ and } ((2k+1)\pi i)^n + \overline{\mu}_j((2k+1)\pi i)^{n-2}$$

respectively. Instead of (14) using

$$(\rho_{k,j} - ((2p+1)\pi i)^n + \mu_s((2k+1)\pi i)^{n-2})\left(X_{k,j}, \widetilde{E}_{p,s}\right)$$
$$= \left((P_2 - C)X_{k,j}^{(n-2)}, \widetilde{E}_{p,s}\right) + \sum_{\nu=3}^{n} \left(P_\nu X_{k,j}^{(n-\nu)}, \widetilde{E}_{p,s}\right)$$

and instead of (60) taking a family of operators $A_\varepsilon = T + \varepsilon(A - T)$, $0 \leq \varepsilon \leq 1$, where $T$ is the operator generated by the differential expression

$$y^{(n)} + (C + (b_{j,j,2k+1}e^{i2\pi(2k+1)x} + b_{j,j,-2k-1}e^{-i2\pi(2k+1)x})I)y^{(n-2)}$$

and boundary conditions (9), and arguing as in the proof of Theorem 3, we get

**Theorem 4** *Let the eigenvalues $\mu_1, \mu_2, ..., \mu_m$ of the matrix $C$ be simple. If the conditions*

$$\lim_{k \to \infty} \frac{\ln|k|}{kb_{s,s,\pm(2k+1)}} = \lim_{k \to \infty} \frac{b_{2k+1}b_{-2k-1}}{b_{s,s,\pm(2k+1)}} = 0, \ \forall s$$

*hold, then:*

*(a) There exists a constant $N_3$ such that the eigenvalues $\rho_{k,j}$ for $|k| \geq N_3$ are simple and consist of 2 sequences $\{\rho_{k,j} : k \geq N_3\}$ and $\{\rho_{-k,j} : k \geq N_3\}$ satisfying*

$$\rho_{\pm k,j} = ((2k+1)\pi i)^n + \mu_j((2k+1)\pi i)^{n-2} \pm ((2k+1)\pi i)^{n-2}q_{2k+1} + O\left(k^{n-3}\gamma_{2k+1}\ln|k|\right).$$

*The corresponding normalized eigenfunction $X_{\pm k,j}(x)$ satisfies*

$$X_{\pm k,j} = (1 + |\alpha_{(2k+1),j}|^2)^{-\frac{1}{2}}(E_{k,j} + \alpha_{\pm(2k+1),j}E_{-k,j}) + O(b_{2k+1}^2) + O(b_{-2k-1}^2) + O\left(\frac{1}{k}\right).$$

*(b) The root functions of $A$ form a Riesz basis in $L_2^m[0,1]$ if and only if $b_{j,j,2k+1} \sim b_{j,j,-(2k+1)}$ for all $j = 1, 2, ..., m$.*

**Acknowledgement 1** *The work was supported by the Scientific and Technological Research Council of Turkey (Tübitak, project No. 108T683).*